\documentclass[reqno]{amsart}
\usepackage{amssymb}
\usepackage{amscd}
\usepackage{mathrsfs,amsmath}
\usepackage{enumitem}
\usepackage{xcolor}

\newtheorem{theorem}{Theorem}[section]

\newtheorem{lemma}{Lemma}[section]

\newtheorem{remark}{Remark}[section]

\numberwithin{equation}{section}

\newcommand{\bigslant}[2]{{\raisebox{.2em}{$#1$}\left/\raisebox{-.2em}{$#2$}\right.}}

\usepackage{xcolor}

\begin{document}
\colorlet{BLUE}{blue}

\title[hyperbolic problem with mixed B.C.]{Initial-boundary value problem for second order hyperbolic operator with mixed boundary conditions}



\author[D. AIT-AKLI]{\bfseries Djamel AIT-AKLI$^{*}$\\\today}

\address{$^1$L2CSP laboratory, Mouloud Mammeri University,  Tizi-Ouzou, 15000, Algeria.}
\address{$^*$Corresponding author: \textnormal{djamel.aitakli@ummto.dz}}

\begin{abstract}
We deal with a linear hyperbolic differential operator of the second order on a bounded planar domain with a smooth boundary. We establish a well-posedness result in case where a mixed, Dirichlet-Neumann, condition is prescribed on the boundary. We focus on the case of a non-homogeneous Dirichlet data and a homogeneous Neumann one. The presented proof is based on a functional theoretical approach and on an approximation argument. Moreover, this work discuss an improvement of a result concerning the range of some operators related to the considered hyperbolic PDE yielding characterizations for the range space of these operators.
\end{abstract}

\maketitle



\section{Introduction}\label{sec1}

Throughout this paper, we let $\Omega\subset\mathbb{R}^2$ to be a bounded planar domain with a $C^\infty$ boundary. This boundary, $\partial \Omega := \Gamma_d \cup \Gamma_n $, is formed out of two parts, each of which is of strictly positive measure. We assume, for simplicity, that $\Gamma_d$ is connected and denote:
\begin{equation}\label{dn}
 \lbrace s_i   \rbrace_{i=\overline{1,2}}  := \overline{\Gamma_d}\cap\overline{\Gamma_n}.
\end{equation}
 We are interested in the study of the second order hyperbolic operator: 
\begin{equation}\label{hyperbol}
\mathcal{H}u := \partial^2_tu - \mathcal{L}u,
\end{equation}
where $\mathcal{L}$ is a second order uniformly and strongly elliptic operator, one can think of the Laplace or the Lam\'e operator for instance. In the sequel, and only for the sake of clarity of presentation, the operator $\mathcal{L}$ will be represented by the Laplace operator $\mathcal{L} \equiv\Delta$. We consider the partial differential equation: 
\begin{equation}\label{eq1}
\mathcal{H}u = f\quad \text{a.e.}~ \text{in}~  (0,T)\times\Omega,
\end{equation}
where $T>0$. The equation (\ref{eq1}) is completed with a mixed boundary condition. More precisely, we prescribe a Dirichlet condition on $(0,T)\times\Gamma_d$ and a Neumann condition on $(0,T)\times\Gamma_n$ i.e. we have in the sens of trace:
\begin{equation}\label{Boundary}
u = G\quad\text{on}~(0,T)\times\Gamma_d\quad,\quad    \frac{\partial u}{\partial \overrightarrow{n}} = 0\quad\text{on}~(0,T)\times\Gamma_n,
\end{equation}
where $\overrightarrow{n}$ is the exterior unit normal defined at each point of $\Gamma_n$ and $\frac{\partial u}{\partial \overrightarrow{n}}$ denotes the normal derivative of $u$. Moreover, initial-in-time conditions are also prescribed: 
 \begin{equation}\label{initial}
  u(0,x) = \Psi^0(x), \quad  \partial_t u(0,x)= \Psi^1(x)\quad \text{on}~\Omega. 
 \end{equation} 
 Such a problem is known to model wave propagation in case when $\mathcal{L} = \Delta$ is the Laplace and models for instance the time-dependent elasticity behavior in case where $\mathcal{L}$ is the Lam\'e system.
We address in this work the issue of existence and uniqueness of a solution to the problem associated to equation (\ref{eq1}) when endowed with the conditions (\ref{Boundary}) and (\ref{initial}). Moreover, we establish an energy estimate that asserts the continuous dependence of the solution of this problem with respect to the data, this estimate expresses the stability of the solution.

Let's place our work in perspective:   a well-posedness result, in case of an homogeneous Dirichlet conditions, can be found in 
\cite{wilcox}. A similar study for the case of a non-homogeneous Dirichlet condition, on the entire boundary, is dealt with in \cite{Lions}. Regarding other types of boundary condition, the work presented in \cite{Trij} deals with the case of the non-homogenous Neumann-type boundary condition. The problem with a mixed homogeneous type condition has been dealt with in \cite{Ibuki}, see also \cite{Hayashida}.

\noindent{\bf Statement of the main result.} The problem under consideration is given in the following system:
\begin{equation}\label{mainprob}
\left\{
    \begin{array}{ll}
\begin{aligned}
&\partial^2_t u - \Delta u = f\quad\text{in}~(0,T)\times\Omega, \\
&            u = G \quad\text{on}~ (0,T)\times\Gamma_d,\\  
&            \frac{\partial u}{\partial \overrightarrow{n}} = 0 \quad\text{on}~ (0,T)\times\Gamma_n, \\
&           u(0,x) = \Psi^0, \quad  \partial_tu(0,x) = \Psi^1~ \quad\text{in}~ \Omega.
\end{aligned}
   \end{array}
\right.
\end{equation}
Let $K(\Omega)$ denote the completion of the set of the up to boundary smooth functions in $\Omega$, vanishing on a neighborhood of $\Gamma_d$, with respect to the Sobolev $H^1(\Omega)$-norm i.e.
\begin{equation}\label{espaceK}
K(\Omega) := \overline{\lbrace v\in C^\infty(\overline{\Omega}), \quad {\rm supp}~v\cap\overline{\Gamma}_d = \emptyset   \rbrace}^{||~||_{H^1}}.
\end{equation}
 Two functions $\xi $ and $ \eta$ are said to satisfy the weak boundary condition, cf. {\cite[Definition 2, p.170]{Ibuki}}, if:
 \begin{equation}\label{wb}
 \begin{aligned}
 &(\xi, \eta)\in  K(\Omega)\times  K(\Omega),\\
 &\int_{\Omega}\Delta \xi\varphi = \int_{\Omega}\nabla \xi\nabla\varphi,\quad      \forall\varphi\in K(\Omega).
 \end{aligned}
   \end{equation} 
   The main result of the paper is stated in the following theorem:
\begin{theorem}\label{mainth}
 Let $f\in C^1([0,T]; L^2(\Omega))$. Assume that $(\Psi^0, \Psi^1)$ satisfy the weak boundary condition (\ref{wb}). Let $G\in H^1((0,T)\times\Gamma_d) $ such that there exists $\tilde{G}\in H^1((0,T)\times\partial\Omega)$ satisfying:
\begin{align}
 &\tilde{G} \equiv G\quad on~ (0,T)\times\Gamma_d,\label{cg1} \\
 &\tilde{G}(0,x) = 0 \quad on ~ \partial\Omega\label{cg2}.
 \end{align} 
  Let $\lbrace s_i\rbrace_i$ be defined as in (\ref{dn}). We assume further that $G$ satisfy point-wisely a domination condition near $s_i$, that is:
 \begin{equation}\label{small}
\exists\gamma >0,~ \lvert G(t, x)\rvert = \underset{x\rightarrow s_i}{\mathcal{O}}(\lvert x - s_i\rvert^{1 + \gamma}), ~\forall t\in (0,T),~ i=1,2.
 \end{equation} 
  Then, there is a unique solution $u\in C^0([0,T];  H^1(\Omega))\cap C^1([0,T]; L^2(\Omega))$ of system (\ref{mainprob}), moreover this solution satisfies the following energy estimate:
\begin{equation}\label{propestim22}
  \begin{split}
&\lVert u\rVert _{C^0([0,T]; H^\alpha(\Omega))} + \lVert\partial_t u\rVert_{C^0([0,T]; H^{\alpha - 1}(\Omega))}\\& \leq C \left( \lVert G\rVert_{H^1((0,T)\times\Gamma_d)}   + \lVert f\rVert_{L^1(0,T; L^2(\Omega))} + \lVert\Psi^0\rVert_{H^1(\Omega)} + \lVert\Psi^1\rVert_{L^2(\Omega)}       \right).
  \end{split}
\end{equation}
\begin{equation}\label{alph}
\text{where}~ \alpha = \frac{3}{5}-\epsilon~ \text{and where} ~\epsilon >0 ~ \text{is arbitrary small}.
\end{equation}

\end{theorem}

\noindent{\bf The plan of the paper.} We prove in section 2 an auxiliary Lemma which states that the solution of a problem similar to (\ref{mainprob}) i.e. with homogeneous second member and initial data, satisfy an energy estimate, cf. Estimate (\ref{propestim22}). As a second step, we construct in section 3 a sequence of functions that solves approximating problems similar to (\ref{mainprob}) and that effectively satisfy energy estimate of Lemma \ref{prop}. We finally pass to the limit and conclude the results of Theorem \ref{mainth}. \\

In all the sequel, we denote by $E^\ast$ the topological dual space of the Banach or Hilbert space $E$ and if $\mathcal{M}$ is a subset of the euclidean space, we denote by $C^\infty_0(\mathcal{M})$ the set of smooth function with compact support in $\mathcal{M}$.

 \section{Preliminary facts} \label{sec:2}
 Before presenting the two auxiliary results, we need to set some preliminary facts: 
 \begin{remark}\label{rem0}
Consider the Dirichlet boundary value problem:
\begin{equation}\label{dirich01}
\left\{
    \begin{array}{ll}
\begin{aligned}
&\partial^2_t u - \Delta u = 0 \quad {\rm in}~ (0,T)\times\Omega, \\
&            u = G \quad {\rm on}~ (0,T)\times\partial\Omega,\\      
&           u(0,x) = 0~ , \quad  \partial_tu(0,x) = 0 \quad {\rm in}~ \Omega,
\end{aligned}
   \end{array}
\right.
\end{equation}   
where $G\in H^1((0,T)\times\partial\Omega)$ satisfy the same compatibility conditions as $\tilde{G}$ in (\ref{cg2}). Following {\cite[Theorem 2.1, p.151]{Lions}}, there exists a unique solution $u \in C^0([0, T]; H^1(\Omega))\cap C^1([0,T]; L^2(\Omega))$ to problem (\ref{dirich01}), moreover $\frac{\partial u}{\partial\overrightarrow{n}}\in L^2((0,T)\times\partial\Omega)$. 
 Let us denote $g :=  \frac{\partial u}{\partial \overrightarrow{n}}\in L^2((0,T)\times\partial\Omega)$ the Neumann data associated to system (\ref{dirich01}). Combining {\cite[Theorem A.1, p. 117]{Trij}} and the remark at the beginning of {\cite[section 1.2, p. 116]{Trij}}, we conclude that $u$ satisfy following energy estimate:
\begin{equation}\label{energy}
\begin{aligned}
&\lVert u \rVert_{C^0\left(0,T; H^\alpha(\Omega)\right)}  + \lVert\partial_t u \rVert_{C^0\left(0,T; H^{\alpha - 1}(\Omega)\right)}\\& \leq C\lVert g \rVert_{L^2((0,T)\times\partial\Omega)}.
\end{aligned}
\end{equation}
where $\alpha$ is given by (\ref{alph}). Using the continuity of the trace operator on $\partial\Omega$ we deduce from (\ref{energy}):
\begin{equation}\label{energyT}
\begin{aligned}
\lVert G\rVert_{C^0\left(0,T; H^{\alpha - \frac{1}{2}}(\partial\Omega)\right)}  = \lVert u\rVert_{C^0\left(0,T; H^{\alpha - \frac{1}{2}}(\partial\Omega)\right)} & \leq C\lVert g\rVert_{L^2((0,T)\times\partial\Omega)}.
\end{aligned}
\end{equation}
\end{remark} 
We state a useful density fact that will be used in the proof of Lemma \ref{prop}:
\begin{remark}\label{dxensi}
Pose  
\begin{equation}\label{J}
\mathcal{J} := \lbrace v\in C^\infty((0,T)\times\partial\Omega),~{\rm supp~}v\cap (\lbrace  0 \rbrace \times\partial\Omega) = \emptyset \rbrace. 
\end{equation}
Pose 
\begin{equation}\label{JJ1}
\mathcal{J}_{\frac{1}{2}} := \overline{\mathcal{J}}^{H^{\frac{1}{2}}}, \quad \mathcal{J}_{1} := \overline{\mathcal{J}}^{H^1},
\end{equation}
 which are respectively the completion of $\mathcal{J}$ with respect to the $H^{\frac{1}{2}}((0,T)\times\partial\Omega)-$norm and $H^1((0,T)\times\partial\Omega)-$norm. We claim that $\mathcal{J}\underset{d}{\hookrightarrow } [\mathcal{J}_{\frac{1}{2}}]^\ast$.  Consider the mixed boundary value problem associated to the Laplace operator on $(0,T)\times\Omega$:
\begin{equation*}
\left\{
    \begin{array}{ll}
\begin{aligned}
&-\Delta v = 0 \quad{\rm in}~ (0,T)\times\Omega, \\ 
&         v = 0 \quad{\rm on}~ \lbrace 0\rbrace\times\Omega, \\
&         \frac{\partial v}{\partial\nu} = g \quad{\rm on}~  ((0,T)\times\partial\Omega) \cup(\lbrace T\rbrace\times\Omega),
\end{aligned}
   \end{array}
\right.
\end{equation*}
where $\partial((0,T)\times\Omega) := ((0,T)\times\partial\Omega) \cup (\lbrace 0\rbrace\times\Omega)\cup (\lbrace T \rbrace\times\Omega)$. Then consider the operator
\begin{alignat*}{2}
K: [\mathcal{J}_{\frac{1}{2}}]^\ast &\rightarrow \mathcal{V} \subset H^1((0,T)\times\Omega) \\
                         g &\mapsto K(g) = v,
\end{alignat*}
where $\mathcal{V} = K([\mathcal{J}_{\frac{1}{2}}]^\ast)$. The operator $K$ is invertible and continuous between the closed, and then Banach, spaces $[\mathcal{J}_{\frac{1}{2}}]^\ast $ and $\mathcal{V}$, we deduce that its inverse $K^{-1}$ is also continuous i.e.
\begin{equation}\label{contJJ}
\lVert K^{-1}(v)\rVert_{[\mathcal{J}_{\frac{1}{2}}]^\ast } =\lVert g\rVert_{[\mathcal{J}_{\frac{1}{2}}]^\ast } \leq C_{-1} \lVert v\rVert_{H^1((0,T)\times\Omega)} 
\end{equation}
 for all $v\in \mathcal{V}$. On the other hand we know that $$\lbrace v\in C^\infty(\overline{(0,T)\times\Omega}),~{\rm supp~}v\cap (\lbrace  0 \rbrace \times\Omega) = \emptyset \rbrace\subset\mathcal{J}$$ is dense in  $\left(\mathcal{V}, \lVert~\rVert_{H^1((0,T)\times\Omega)}\right)$ and that $K^{-1}$ maps $C^\infty$ functions into $C^\infty$ functions. Thus we deuce, using (\ref{contJJ}), that $\mathcal{J}$ is dense in $[\mathcal{J}_{\frac{1}{2}}]^\ast$. 
Let us show that 
\begin{equation}\label{densC1}
C^\infty_0((0,T)\times\partial\Omega)\underset{d}{\hookrightarrow} [H^1((0,T)\times\partial\Omega)]^\ast.
\end{equation}
 Recall the standard result $C^\infty_0((0,T)\times\partial\Omega)\underset{d}{\hookrightarrow}L^2((0,T)\times\partial\Omega)$. On another hand, for any element $f\in [H^1((0,T)\times\partial\Omega)]^\ast$ there exists $f_0, f_1, f_2, f_3\in L^2((0,T)\times\partial\Omega)$ such that $f$ is represented by: $ <f,v> = \int f_0 v + f_1v_t + f_2{x_1} + f_3 v_{x_2}.$ The density claim follows immediately.
To obtain such a representation of $f\in [H^1]^\ast$, one can use the Riesz representation theorem and the definition of the Sobolev space $H^1(\mathcal{M})$ on the manifold $\mathcal{M}:= (0,T)\times\partial\Omega$, see for instance {\cite[Section 2.2, p.21-22]{manifsobo}}.
\end{remark}
 \section{Auxiliary results} \label{sec:3}
 We state the essential auxiliary Lemma needed subsequently for the proof of Theorem \ref{mainth}:
\begin{lemma}\label{prop}
 Consider the mixed boundary value problem given by:
\begin{equation}\label{mix}
\left\{
    \begin{array}{ll}
\begin{aligned}
&\partial^2_t u - \Delta u = 0 \quad\text{in}~ (0,T)\times\Omega, \\
&            u = G  \quad\text{on}~ (0,T)\times\Gamma_d,\\  
&  \frac{\partial u}{\partial \overrightarrow{n}} = 0  \quad\text{on}~ (0,T)\times\Gamma_n, \\
&           u(0,x) = 0, \quad  \partial_tu(0,x) = 0  \quad\text{in}~ \Omega. 
\end{aligned}
   \end{array}
\right.
\end{equation}
where we assume that $G$ meet the same compatibility conditions as $\tilde{G}$ in (\ref{cg2}) on $\Gamma_d$. Then the subspace $\mathcal{G} \subset H^1((0,T)\times\Gamma_d)$: 
\begin{equation}\label{subG}
\mathcal{G} :=
\left\{\!\begin{aligned}
&G\in H^1((0,T)\times\Gamma_d): \text{problem}~(\ref{mix})~\\& \text{admits a}\text{ unique solution} ~u~such~ that ~\\ & u\in C^0(0,T; H^\alpha(\Omega))\cap C^1(0,T; H^{\alpha - 1}(\Omega))~
\end{aligned}\right\}.
\end{equation}
isn't trivial. Moreover, any solution $u$ of problem (\ref{mix}) satisfy the energy estimate:
\begin{alignat}{2}\label{propestim}
&\lVert\frac{\partial u}{\partial \overrightarrow{n}}\rVert_{H^\beta((0,T)\times \Gamma_d)} + \lVert u\rVert_{C^0(0,T; H^\alpha(\Omega))} + \lVert\partial_t u\rVert_{C^1(0,T; H^{\alpha - 1}(\Omega))}\nonumber\\& \leq C_1  \lVert G\rVert_{H^1((0,T)\times\Gamma_d)},
\end{alignat}
for all $G\in\mathcal{G}$, where the exponent $\beta$ equals $2(1-\alpha)$ and $\alpha$ is given by (\ref{alph}). A trivial subspace means: a vector space containing only the zero element  $\mathbf{0}_{H^1}$. 
 \end{lemma}

\noindent Regarding the non triviality of $\mathcal{G}$, we construct, in step 1 of the proof of the main result cf. problem (\ref{sys11}), a sequence of functions that are different from the zero element and which belong to $\mathcal{G}$, this shows the first claim of Lemma \ref{prop}. Remark that:
\begin{equation}\label{D0}
\left\{\!\begin{aligned}
&G\in H^1((0,T)\times\partial\Omega): G ~ \text{ satisfy the} \\ &\text{  compatibility conditions}~ (\ref{cg2})~\text{on}~\partial\Omega\end{aligned}\right\}\equiv \overline{\mathcal{J}}^{H^1}.
\end{equation}
where $\mathcal{J}$ is defined by (\ref{J}). The equivalence is due to the fact that $\partial\left((0,T)\times\partial\Omega \right) = \lbrace\lbrace 0\rbrace \times\partial\Omega\rbrace\cup \lbrace\lbrace 1\rbrace \times\partial\Omega \rbrace$. Consider the non necessarily homogeneous Dirichlet problem given by the following system:
\begin{equation}\label{dirich}
\left\{
    \begin{array}{ll}
\begin{aligned}
&\partial^2_t u - \Delta u = 0\quad\text{in}~ (0,T)\times\Omega, \\
&            u = G \quad\text{on}~ (0,T)\times\partial\Omega,\\      
&           u(0,x) = 0~,~  \partial_tu(0,x) = 0\quad\text{in}~ \Omega,
\end{aligned}
   \end{array}
\right.
\end{equation}
where $ G  \in\mathcal{J}_1\subset H^1((0,T)\times\partial\Omega)$. Consider the operator $T$ defined by:
\begin{equation}\label{TT}
\begin{aligned}
T: H^1((0,T)\times\partial\Omega)&\rightarrow  L^2((0,T)\times\partial\Omega)\\
                                 G&\mapsto T(G) = g:=  \frac{\partial u}{\partial \overrightarrow{n}},
\end{aligned}
\end{equation}
the operator $T$ associates to every Dirichlet data $G\in \mathcal{J}_1\subset H^1((0,T)\times\partial\Omega)$, the unique Neumann data $g \in  L^2((0,T)\times\partial\Omega)$ corresponding to the solution $u$ of problem (\ref{dirich}). Following {\cite[Theorem 2.1, p.151]{Lions}}, the operator $T$ is well defined. According to {\cite[Remark 2.2, p.152]{Lions}} the operator $T: H^1\rightarrow L^2$ is also bounded. That said, $T(H^1((0,T)\times\partial\Omega))\subset L^2((0,T)\times\partial\Omega)$ can not be complete for the $L^2$-norm. This is the reason why we restrict, in the sequel, the operator $T$ to more regular spaces, namely $\mathcal{D}$ and $\mathcal{A}$, which are respectively given in (\ref{D}) and (\ref{A0}). But before doing this, we need to consider less regular space. For the case of distributional Neumann data, following {\cite[Theorem G, p.119]{Trij}}, we know that 
\begin{equation}\label{eqin}
   T^{-1}([H^1((0,T)\times\partial\Omega)]^\ast) \subset [Y_\beta]^\ast,
\end{equation}
where
\begin{equation}\label{Ybeta}
 Y_\beta := \lbrace h\in H^\beta((0,T)\times\partial\Omega),~ h(0,.) = 0 \rbrace,
\end{equation} 
 and  $\beta = 2(1-\alpha)$, where $\alpha$ is given by (\ref{alph}).
  
 Before going further, we need to derive a crucial adjoint identity, this is the subject of Lemma \ref{adjo}:
\begin{lemma}\label{adjo}
Let $Y_\beta$ be given by (\ref{Ybeta}). Two functions $G \in T^{-1}(T(H^1)\cap Y_\beta)\subset H^1((0,T)\times\partial\Omega)$ and $g\in T(H^1)\cap Y_\beta$ are respectively the Dirichlet and Neumann data relatively to a solution of problem (\ref{dirich}) if and only if the following adjoint identity is fulfilled:
\begin{equation}\label{adjoF}
< G, \eta>_{1,[1]^\ast}   = < g,  h>_{Y_\beta, [Y_\beta ]^\ast} 
\end{equation}
for all $\eta \in [H^{1}((0,T)\times\partial\Omega)]^\ast$ and all $h = T^{-1}(\eta)\in [Y_{\beta}]^\ast $.
\end{lemma}
 We will prove in the sequel, namely as a consequence of the first step of the proof of Lemma \ref{prop}, that  $T(\mathcal{J}_1)\equiv Y_\beta$, thus estimate (\ref{adjoF}) holds for $G\in \mathcal{J}_1$ and $g\in T(\mathcal{J}_1) = Y_\beta$, where $\mathcal{J}$ and $\mathcal{J}_1$ are given by (\ref{JJ1}). We now present a proof of Lemma \ref{adjo}:
\begin{proof}
Consider the problem, cf. {\cite[Problem 2.63, p.163]{Lions}}:
\begin{equation}\label{phi}
\left\{
    \begin{array}{ll}
\begin{aligned}
&\partial^2_t \phi - \Delta \phi = 0\quad\text{in}~ (0,T)\times\Omega, \\
&            \phi = h \quad\text{on}~ (0,T)\times\partial\Omega,\\      
&          \phi(T,x) = \partial_t \phi(T,x) = 0\quad\text{in}~ \Omega.
\end{aligned}
   \end{array}
\right.
\end{equation}
Pose: $\eta(t,x) := \frac{\partial\phi}{\partial\overrightarrow{n}}(t,x) $. According to {\cite[Estimate (2.65), p.163]{Lions}}, we know that if $G\in\mathcal{J}_1 \subset H^1((0,T)\times\Omega)$ and $g\in L^2$ are respectively the Dirichlet and Neumann data corresponding to the  solution to problem (\ref{dirich}), then we have:
\begin{equation}\label{use}
\int_{(0, T)\times\partial\Omega} G \eta   = \int_{(0, T)\times\partial\Omega}g h. 
\end{equation}
 for all $ h\in H^1((0,T)\times\partial\Omega)\supset H^1_0((0,T)\times\partial\Omega)$  and  all $\eta = T( h ) \in  L^2((0,T)\times\partial\Omega)$ such that $(h, \eta)$ corresponds to the solution $\phi$ of problem (\ref{phi}), which presupposes the compatibility condition $h(T,.) = 0$, cf. Problem (\ref{phi}). 
   Set $\mathcal{F}:= \lbrace \phi\in H^1((0,T)\times\Omega): \phi(T,.) = \partial_t \phi(T,.) = 0 \rbrace $. We claim that: $C^\infty_0((0,T)\times\partial\Omega) \subset T(\mathcal{F}\restriction_{(0,T)\times\partial\Omega)})$ i.e. if the Neumann data is $C^\infty_0$, then the corresponding solution $u$ of the Dirichlet problem (\ref{dirich}) satisfy $u\in\mathcal{F}$. To see it, consider the Neumann problem given by:
\begin{equation}\label{eta}
\left\{
    \begin{array}{ll}
\begin{aligned}
&\partial^2_t \varphi - \Delta \varphi = 0\quad\text{in}~ (0,T)\times\Omega, \\
&            \nabla\varphi\cdot\overrightarrow{n} = \eta \quad\text{on}~ (0,T)\times\partial\Omega,\\      
&          \varphi(0,x) = \partial_t \varphi(0,x) = 0\quad\text{in}~ \Omega.
\end{aligned}
   \end{array}
\right.
\end{equation}
We denote, in case when $\phi$ is smooth, $h := \varphi\restriction_{(0,T)\times\partial\Omega}$. There is an important fact to emphasize. Let $\eta \in C^\infty_0((0,T)\times\partial\Omega)$. Consider the $T-$shifted function $t\mapsto\eta_T(t) := \eta(T-t)$. Then $\eta(0) = \eta_T(0) = 0$. Since the given Neumann data $\eta$ are mapped into solutions $\varphi$ of problem (\ref{eta}), then $T-$shifted Neumann data $\eta(T-.,.)$ are mapped into $T-$shifted solution $\varphi_T$; to see this, one may consider {\cite[Identity 1.3, 1.6 p.115]{Trij}}. Consequently for every $\eta\in C^\infty_0((0,T)\times\partial\Omega)$, the solution $\varphi_T$ corresponding to $\eta_T$ satisfy $ \varphi_T \in\mathcal{F}$. But since $\eta \in C^\infty_0$ if and only if $\eta_T \in C^\infty_0$, then for every $\eta\in C^\infty_0((0,T)\times\partial\Omega)$, the corresponding solution $\varphi$  satisfy $\varphi \in \mathcal{F}$, this proves the claim. We infer that every $(h, \eta)$ satisfying $\eta\in C^\infty_0((0,T)\times\partial\Omega)$ and $h = T^{-1}(\eta)$ fulfills the identity (\ref{use}).  
 On another hand, following {\cite[Theorem G, p.119]{Trij}}, one can see that for every Neumann data $\eta\in [H^1\left((0,T)\times\partial\Omega\right)]^\ast$, there exists a unique solution to problem (\ref{eta}) whose corresponding Dirichlet data satisfy  $h = T^{-1}(\eta)\in [Y_{\beta}]^\ast$. Moreover, one may consider the continuity statement at the beginning of {\cite[section 1.2, p. 116]{Trij}} to deduce that there exists $C>0$ s.t.:
\begin{equation}\label{etaC}
\lVert h\rVert_{[Y_\beta]^\ast((0,T)\times\partial\Omega)} = \lVert T^{-1}(\eta)\rVert_{[Y_\beta]^\ast}\leq C \lVert\eta\rVert_{[H^{1}\left((0,T)\times\partial\Omega\right)]^\ast},
\end{equation}
for all $\eta\in [H^1\left((0,T)\times\partial\Omega\right)]^\ast$ where $Y_\beta$ is given by (\ref{Ybeta}). According to the claim proved in the beginning of this proof, estimate (\ref{use}) can be rewritten in the following setting: for $G \in T^{-1}\left( T(H^1)\cap Y_\beta\right)\subset H^1$ and $g\in Y_\beta\cap T(H^1)$ we have: 
\begin{equation}\label{use1}
< G, \eta>_{1,[1]^\ast}   = < g,  h>_{Y_\beta, [Y_\beta]^\ast} 
\end{equation}
for all $\eta\in C_0^\infty((0,T)\times\partial\Omega)$ and all $h = T^{-1}(\eta)$, $h(T,.) = h(0,.) = 0$. Let $\eta\in [H^1((0,T)\times\partial\Omega)]^\ast$. Applying the density fact (\ref{densC1}), stated in Remark \ref{dxensi}, there exists a sequence $(\eta^m)_m $ in $C^\infty_0((0,T)\times\partial\Omega)$ such that $\lVert\eta^m -\eta\rVert_{[H^{1}]^\ast}\rightarrow 0$, this convergence implies, by using (\ref{etaC}), that $\lVert h^m - h\rVert_{[Y_\beta]^\ast}\rightarrow 0$, where $h^m := T^{-1}(\eta^m)$. Then $(h^m, \eta^m)$ do satisfy:
\begin{equation}\label{use15}
< G, \eta^m>_{1,[1]^\ast}   = < g,  h^m>_{Y_\beta, [Y_\beta]^\ast} 
\end{equation}
for all $m\in\mathbb{N}$. Letting $m\rightarrow \infty$ in (\ref{use15}), we discover that (\ref{use}) can be extended to the form:
\begin{equation}\label{use2}
< G, \eta>_{1,[1]^\ast}   = < g,  h>_{Y_\beta, [Y_\beta]^\ast} 
\end{equation}
for all $\eta \in [H^{1}((0,T)\times\partial\Omega)]^\ast$ and all $h = T^{-1}(\eta)\in [Y_{\beta}]^\ast $. The estimate (\ref{use2}) is the desired adjoint identity. 
Actually we have established the necessary condition part of the lemma. To see the other direction of the equivalence, it suffices to integrate by part in (\ref{use2}) against adequately regular test functions $\eta $ and $h$ and then deduce that $(G, g)$ actually solves problem (\ref{dirich}). 
\end{proof}

The identity (\ref{adjoF}) in Lemma \ref{adjo} is the key for proving the completeness of the functional space $\mathcal{D}$ defined by (\ref{D}). Denote $\mathcal{E}$ to be the closure of $ T^{-1}([H^1]^\ast)$ with respect to the strong norm of the dual, $b(Y_\beta, [Y_\beta]^\ast)$:
\begin{equation}\label{E}
\mathcal{E} := \overline{ T^{-1}([H^1]^\ast)}^{b([Y_\beta], [Y_\beta]^\ast)}\subset [Y_\beta]^\ast.
\end{equation}
 Notice that, since the involved spaces in the topology $b$ are normed, this dual topology is equivalent to the topology induced by the operator norm. The space $\mathcal{E}$ is a Banach and admits a (non necessarily unique) pre-dual. The existence of such a pre-dual, that we denote $\mathcal{A}$, can be shown using the sufficient condition stated in {\cite[Proposition 1, p.88]{predual}}. Indeed, since $\mathcal{E}$ is a complete, and so a closed subspace of the reflexive space $[Y_\beta]^\ast$, then it is also reflexive. This means that the natural embedding in the bidual is onto  an thus it can be easily shown that the, above mentioned, sufficient condition for existence is fulfilled.
Let us denote   
\begin{equation}\label{A0}
\mathcal{A} := \text{a pre-dual of the space}~ \mathcal{E}~ \text{containing }~ \mathcal{J}\supset C^\infty_0.
\end{equation}
where $\mathcal{J}$ is defined by (\ref{J}). We report four observations:\\

{\it\paragraph{\bf Observation 1.} We claim that $\mathcal{A}\subset Y_\beta$. Consider the isometry $\Psi:  Y_\beta\rightarrow [Y_\beta]^\ast$ that associate to each element of $Y_\beta$ the unique corresponding element in $[Y_\beta]^\ast$ via the Riesz representation theorem with respect to the inner product induced by $( , )_{H^\beta} $ i.e. for every $f\in [Y_\beta]^\ast$, there exists a unique $u_f\in Y_\beta$ such that:  $(u_f, v )_{H^\beta} = <f, v> $ for all $v\in Y_\beta$. Since $\mathcal{E}$ is a closed subspace of $[Y_\beta]^\ast$, then by continuity of $\Psi$ we deduce that $ \Psi^{-1}(\mathcal{E})$ is also closed and thus is a Hilbert space. By definition, every element of $\mathcal{E}\subset [Y_\beta]^\ast$ is continuous on $Y_\beta$, and thus is continuous on $\Psi^{-1}(\mathcal{E})\subset Y_\beta$, this is due to the fact that both $Y_\beta$ and  $\Psi^{-1}(\mathcal{E})$ are endowed with the same norm $\lVert~\rVert_{H^\beta}$. We consider the restriction of $\Psi$ defined by:
$\Psi:  \Psi^{-1}(\mathcal{E})\rightarrow \mathcal{E}$, it is still an isometric isomorphism. Using (\ref{A0}), we infer that $\mathcal{A}\cong \Psi^{-1}(\mathcal{E})$. Hence $\mathcal{A}\subset Y_\beta$ and $\mathcal{A}$ is closed, thus it a Banach space.\\

}

 Pose:
\begin{equation}\label{D}
\mathcal{D} := \lbrace G\in \mathcal{J}_1, : \text{ there exists}~ g\in \mathcal{A},~~ T(G) = g \rbrace.
\end{equation}
where $\mathcal{J}_1 := \overline{\mathcal{J}}^{H^1}$. The subspace $\mathcal{D}$ is endowed with the $H^1((0,T)\times\partial\Omega)$-norm. This space is characterized as the largest subspace of $\mathcal{J}_1$ whose direct image by $T$ is included in $\mathcal{A}$ i.e. such that $\mathcal{A}\supset T(\mathcal{D})$.

{\it \paragraph{\bf Observation 2.} Let $\mathcal{A}$ be a predual of $\mathcal{E}$. We claim that $\mathcal{J}_1 := \overline{\mathcal{J}}^{H^1}\underset{c}{\hookrightarrow }\mathcal{A}$. Indeed, first recall that the spaces $\mathcal{A}$ and $\mathcal{J}_1$ are reflexive and Banach. Moreover, we have clearly (since $ \mathcal{J}\underset{d}{\hookrightarrow } Y_\beta $ and $\mathcal{J}\underset{c}{\hookrightarrow } \mathcal{J}_{\frac{1}{2}}\underset{c}{\hookrightarrow } Y_\beta$) that  $\mathcal{J}_{\frac{1}{2}}\underset{d}{\hookrightarrow } Y_\beta $, then $[Y_\beta]^\ast \underset{c}{\hookrightarrow } [\mathcal{J}_{\frac{1}{2}}]^\ast$ and thus, following (\ref{E}), we have $\mathcal{E} \underset{c}{\hookrightarrow } [\mathcal{J}_{\frac{1}{2}}]^\ast$. According to the density fact stated in Remark \ref{dxensi}, we have $\mathcal{J} \underset{d}{\hookrightarrow } [\mathcal{J}_{\frac{1}{2}}]^\ast$. Consequently, and since $\mathcal{J} \subset\mathcal{E}$,  we have $\mathcal{E}\underset{d}{\hookrightarrow }[\mathcal{J}_{\frac{1}{2}}]^\ast$. We conclude immediately that $\mathcal{J}_{\frac{1}{2}} \underset{c}{\hookrightarrow } [\mathcal{E}]^\ast\equiv \mathcal{A}.$ Since $ \mathcal{J}_1 \underset{c}{\hookrightarrow }\mathcal{J}_{\frac{1}{2}}$, the claim follows. }
 
 {\it \paragraph{\bf Observation 3.}
 The subspace $\mathcal{A}$, being a predual of the closed and then reflexive subspace $\mathcal{E}\subset [Y_\beta]^\ast$, is reflexive. Moreover, according to observation 1, it is complete (with respect to the norm $||~||_{H^\beta}$). From preceding observations $\mathcal{J}\subset\mathcal{J}_1\subset\mathcal{A}\subset Y_\beta$, then we immediately deduce that:
\begin{equation}\label{identYA}
\mathcal{A} \equiv Y_\beta.
\end{equation}
This implies that $\mathcal{E} := [\mathcal{A}]^\ast = [Y_\beta]^\ast$ then, using (\ref{E}), we have 
\begin{equation}\label{densityB}
T^{-1}([H^1]^\ast)\underset{d} {\hookrightarrow} [Y_\beta]^\ast.
\end{equation}
  We claim that actually we have  $ T^{-1}([H^1]^\ast)^\ast \cong [Y_\beta]^{\ast\ast}\cong [Y_\beta]$ (since $Y_\beta$ is reflexive). Indeed, since $T^{-1}([H^1]^\ast)$ is a dense subspace of the Banach space $[Y_\beta]^\ast $ then its orthogonal complement is trivial i.e. $ T^{-1}([H^1]^\ast)^\bot = \lbrace 0 \rbrace$. Consequently, we have the isometric isomorphism: 
\begin{equation}\label{isofor}
[Y_\beta] \equiv \bigslant{[Y_\beta]}{\left( T^{-1}([H^1]^\ast)\right)^\bot}\cong  [T^{-1}([H^1]^\ast)]^\ast.
\end{equation}
this shows the claim.}
{\it \paragraph{\bf Observation 4.}
According to {\cite[Remark 2.10, p.167]{Lions}},  for every $G\in \mathcal{J}$,  we have $g = T(G) \in C^\infty((0,T)\times\partial\Omega)$. Moreover, given the necessary compatibility condition (of the Neumann data with initial data) we infer that $g(0,x) = \frac{\partial u}{\partial\nu}(0,x) = 0 $ for $x\in\partial\Omega$; 
  indeed, this can be deduced using the continuity of the gradient of the solution corresponding to $G$ on the surface $ (0,T)\times\partial\Omega $ (near $\lbrace  0 \rbrace \times\partial\Omega)$, cf. the classical theory of Neumann problems for linear Hyperbolic operators.
Consequently, $g = T(G)\in \mathcal{A}$ i.e.  $ T(\mathcal{J})\subset \mathcal{A}$, we infer from these arguments that $\mathcal{J}\subset\mathcal{D}$.}\\

 An important consequence of this observation is the existence of a unique predual $\mathcal{A}$ of the space $\mathcal{E}$ containing $\mathcal{J}$. Observation 2 shows also that $\mathcal{D}$ and $\mathcal{A}$ aren't trivial. 

The operator $T:  \mathcal{D}\rightarrow T(\mathcal{D})\subset\mathcal{A}$, given by (\ref{TT}), is clearly linear and invertible. The procedure for proving estimate (\ref{propestim}) is based on the boundedness of its inverse operator, denoted $T^{-1}$. For this we need to make use of the Banach isomorphism theorem. In order to apply this theorem, we need, in addition to the fact that  $(\mathcal{A}, \Vert~ \Vert_{H^\beta})$ is complete, to show that $\mathcal{A} = T(\mathcal{D})$ and that the subspace $(\mathcal{D},  \Vert~ \Vert_{H^1})$ is a Banach spaces. We now present a proof of the second claim of Lemma \ref{prop}:

Recall that Lemma \ref{adjo} asserts that: $G\in \mathcal{D}$ and $g = T(G)\in \mathcal{A}$ are respectively the Dirichlet and Neumann data associated to the solution of problem (\ref{dirich}) if and only if (\ref{adjoF}) holds for all $\eta \in [H^1]^\ast $  and all $h\in T^{-1}([H^1]^\ast)\subset [\mathcal{A}]^\ast\equiv \mathcal{E}$.

\subsubsection*{Proof of Lemma \ref{prop}:}
\begin{proof}
The proof is done in two steps:
\paragraph{\bf Step 1: completeness of $\mathcal{D}$.} First recall that $\mathcal{A} \equiv Y_\beta$ is reflexive. Let $\mathcal{D}$ be as defined by (\ref{D}). Let $G_n\in \mathcal{D}$ be a convergent sequence i.e. such that $\Vert G_n - G\Vert_{H^1}\rightarrow 0$, then we show that $G\in\mathcal{D}$. By assumption, for all $n\in\mathbb{N}$, there exists $g_n\in \mathcal{A}\equiv Y_\beta$ such that $T(G_n) = g_n$. Surely we have the weak convergence:
\begin{equation}\label{ConG}
 G_n\underset{H^1}{\rightharpoonup} G.
\end{equation}
 Using (\ref{adjoF}) we can write:
\begin{equation}\label{cov}
< G_n , \eta   >_{1, [1]^\ast} = <g_n , h>_{Y_\beta, [Y_\beta]^\ast}, 
\end{equation}
for all $\eta\in [H^{1}]^\ast$ and all $h = T^{-1}(\eta)\in T^{-1}([H^{1}]^\ast)$.
Then, combining (\ref{ConG}) and (\ref{cov}), we infer that 
\begin{alignat*}{2}
g_n: T^{-1}([H^{1}]^\ast) &\rightarrow  \mathbb{R}   \\
                        h  &\mapsto  <g_n , h>_{Y_\beta, [Y_\beta]^\ast}
\end{alignat*}
 is point-wisly bounded i.e. for all $h\in  T^{-1}([H^{1}]^\ast) $ there is $c_h>0$ such that 
\begin{equation}\label{bnb}
 \sup_n \vert <g_n , h>_{Y_\beta, [Y_\beta]^\ast}\vert \leq c_h.
\end{equation}  
  According to (\ref{isofor}), $g_n\in [Y_\beta] \cong [T^{-1}([H^{1}]^\ast)]^\ast$, then $g_n$ can be viewed as a linear functional on $T^{-1}([H^{1}]^\ast)$. On another hand, since $T^{-1}([H^{1}]^\ast) \underset{d}{\hookrightarrow}[Y_\beta]^\ast$ then, for all $n\in\mathbb{N}$, there exists a unique bounded linear functional $F_n$ that extends $g_n$ from $ T^{-1}([H^1]^\ast)$  to the Banach space $[Y_\beta]^\ast $ such that $ \Vert F_n\Vert_{\mathcal{L}( [Y_\beta]^\ast ,\mathbb{R})} = \Vert g_n\Vert_{\mathcal{L}( T^{-1}([H^{1}]^\ast) ,\mathbb{R} ) }$ i.e. such that
\begin{equation}\label{egalnorm}
\sup_{h\in T^{-1}([H^{1}]^\ast)}\frac{< g_n , h  >}{\lVert h\rVert_{    [Y_\beta]^\ast   }} =                            \sup_{h\in [Y_\beta]^\ast}\frac{< F_n , h  >}{\lVert h\rVert_{[Y_\beta]^\ast}}.
\end{equation} 
Following (\ref{bnb}) We deduce that $F_n$ is point-wisly bounded on $[Y_\beta]^\ast $. Since $[Y_\beta]^\ast$ is a Banach space then, according to the uniform boundedness principle, we deduce that $F_n $ is uniformly bounded. Thus, using (\ref{isofor}) and (\ref{egalnorm}), we infer that 
  $\Vert g_n\Vert_{[ T^{-1}([H^{1}]^\ast)]^\ast} = \Vert g_n \Vert_{Y_\beta}  $ is uniformly bounded. i.e. there exists $g\in Y_\beta$ such that $g_n\underset{Y_\beta}{\rightharpoonup} g$.
Then, letting $n\rightarrow\infty$ in (\ref{cov}),
we deduce easily that $\mathcal{A} \ni g = T(G)$ and thus $G\in \mathcal{D}$. Consequently $\mathcal{D}$ is complete. Hence, $\left(\mathcal{D}, \Vert ~\Vert_{H^1}\right)$ is a Banach space. Using {\it Observation 4} we conclude, by density, that $\mathcal{D}\equiv \mathcal{J}_1$. As a consequence of the preceding argument, we infer that weakly convergent sequences in $\mathcal{J}_1$ are mapped into weakly convergent sequences in $\mathcal{A}$. This observation leads to the conclusion that the operator $T: \left(\mathcal{J}_1, \vert~\Vert_{H^1} \right)\rightarrow \left(\mathcal{A}, \Vert~\Vert_{H^\beta}   \right)$ is weakly continuous, and since it is linear, then it is bounded i.e. $\exists C>0$:
\begin{equation}\label{bounded}
 \Vert T(G)\Vert_{H^\beta( (0,T)\times\partial\Omega )}   \leq C \Vert G\Vert_{H^1((0,T)\times\partial\Omega)},
\end{equation}
 for all $ G\in \mathcal{J}_1$.
 
 Another important consequence of the above arguments is that $T(\mathcal{J}_1) = Y_\beta$. Indeed, we know that $\mathcal{J}\underset{d}{\hookrightarrow }\mathcal{A}$. Let $g\in\mathcal{A}$, there exists $(g_n)_n\in \mathcal{J}$ such that $\Vert g - g_n\Vert_{H^\beta}\rightarrow 0 $, according to the above results, there exists $(G_n)_n  \in  \mathcal{J}_1$ such that $g_n = T(G_n)$ for all $n$. Following these assumptions, $<g_n, h>_{}$ is a Cauchy sequence for all $h \in  T^{-1}([H^{1}]^\ast)$, thus  $<G_n, \eta >$ is also Cauchy for all $\eta\in [H^1]^\ast$. Then, since $\mathcal{D} \equiv \mathcal{J}_1$ is weakly complete (it is indeed reflexive), there exists $G\in H^1$ such that:
 $$<G, \eta   > = <g , h   >    $$ for all $\eta\in [H^1]^\ast$ and all $h\in T^{-1}([H^1]^\ast)$. We deduce easily, by performing an integration by part, that $g = T(G)$, this implies that $g\in T(\mathcal{D})$, thus  $T(\mathcal{J}_1) = \mathcal{A}$.
 
\paragraph{\bf Step 2: proof of Estimate (\ref{propestim}).}  Since $\mathcal{J}_1$ and $\mathcal{A}$ are Banach space and using (\ref{bounded}), we apply the isomorphism theorem, to infer that the inverse operator is bounded i.e. there exists $C_{-1}>0$ s.t.
\begin{equation}\label{continv}
\Vert T^{-1}(g)\Vert_{H^1((0,T)\times\partial\Omega)} = \Vert G\Vert_{H^1((0,T)\times\partial\Omega)}\leq C_{-1}\Vert g\Vert_{H^\beta((0,T)\times\partial\Omega)},
\end{equation}
for all $g\in T(\mathcal{J}_1)\equiv \mathcal{A}\equiv Y_\beta$. Set $\mathcal{B}$ to be the subspace: 
$$ \mathcal{B} := \lbrace  g\in H^\beta((0,T)\times\partial\Omega): g \equiv 0 ~\text{on}~(0,T)\times \Gamma_n    \rbrace.$$
The subspace $\left(\mathcal{B}, \Vert ~\Vert_{H^\beta}\right)$ is clearly closed. We infer that  $\mathcal{A}\cap \mathcal{B}$ is also closed, thus it is a Banach space. Since $T: \mathcal{J}_1\subset H^1 \rightarrow \mathcal{A}$ is continuous, as well as its inverse, then $T^{-1}(\mathcal{A}\cap \mathcal{B})$ is a complete subspace of $\mathcal{J}_1$, thus it is a Banach. We remark that:

\begin{itemize}
\item  The space $ T^{-1}(\mathcal{A}\cap \mathcal{B})$ isn't trivial. Indeed, one can see with ease that the non trivial space $\mathcal{G}$, given by (\ref{subG}), is nothing but the restriction to $\Gamma_d$ of functions of the spaces $T^{-1}(\mathcal{A}\cap \mathcal{B})$, this is due to the fact that any solution $u$ of problem (\ref{mix}) should satisfy the compatibility conditions (\ref{cg2}) on $(0,T)\times\partial\Omega$, and thus we infer the claimed non triviality.
\end{itemize}

\noindent We immediately deduce from estimate (\ref{continv}) that: 
\begin{equation}\label{estiss}
\begin{aligned}
\Vert G\Vert_{H^1((0,T)\times\Gamma_d)}& = \Vert T^{-1}(g)\Vert_{H^1((0,T)\times\Gamma_d)}\\ & \leq C_{-1}\Vert g\Vert_{H^\beta((0,T)\times\Gamma_d)},
\end{aligned}
\end{equation}
for all $g\in\mathcal{A}\cap\mathcal{B}$. Then the restriction of the operator $T^{-1}$ defined by:
\begin{alignat*}{2}
T^{-1}: \mathcal{A}\cap \mathcal{B} &\rightarrow T^{-1}(\mathcal{A}\cap \mathcal{B})\subset \mathcal{J}_1\\
              g &\mapsto T^{-1}(g) = G
\end{alignat*}
 is continuous. Consider the operator ${\rm Re}_{\lvert(0,T)\times\Gamma_d}$ defined by:
\begin{alignat*}{2}
{\rm Re}_{\lvert(0,T)\times\Gamma_d}:  T^{-1}(\mathcal{A}\cap \mathcal{B})\subset \mathcal{J}_1 &\rightarrow \mathcal{C}\\
           G &\mapsto {\rm Tr}_{\lvert(0,T)\times\Gamma_d}G,
\end{alignat*}
where $\mathcal{C} := {\rm Re}_{\lvert(0,T)\times\Gamma_d}\circ T^{-1}(\mathcal{A}\cap\mathcal{B})$. The space $\mathcal{C}$ is the set of restrictions to $(0,T)\times\Gamma_d$ of Dirichlet data $G\in  T^{-1}(\mathcal{A}\cap \mathcal{B})$ that correspond to Neumann data, in $Y_\beta$, which are identically zero on $\Gamma_n$. Following the above remark, we have:
 $${\rm Re}_{\lvert(0,T)\times\Gamma_d}\circ T^{-1}(\mathcal{A}\cap \mathcal{B})\equiv\mathcal{G},$$ where $\mathcal{G}$ is given by (\ref{subG}).
 Estimate (\ref{estiss}) shows that the operator $$ P := \left( {\rm Re}_{\lvert(0,T)\times\Gamma_d}\circ T^{-1}\right) : \mathcal{A}\cap\mathcal{B}\rightarrow\mathcal{C} \equiv \mathcal{G}$$ is continuous, and it is, obviously, one-to-one and onto. One can easily see, using the completeness of  $\mathcal{A}$ and $\mathcal{B}$, in addition to the estimate (\ref{bounded}), that the subspace $\mathcal{C}\equiv \mathcal{G}$ is a Banach space. Then applying the isomorphism theorem we show that the operator $P$ has a bounded inverse  i.e. there exists $C_1>0$ s.t.:
\begin{equation}\label{ess}
 \Vert g\Vert _{H^\beta((0,T)\times\Gamma_d)} = \Vert P^{-1}(G) \Vert_{L^2((0,T)\times\Gamma_d)} \leq C_1\Vert G\Vert_{H^1((0,T)\times\Gamma_d)},
\end{equation}
for all $G\in \mathcal{G} \subset H^1((0,T)\times\Gamma_d) $. Finally, by combining estimates (\ref{ess}), (\ref{energy}) and the continuous mebedding of $H^\beta$ in $L^2$, we conclude immediately Estimate (\ref{propestim}). This ends the proof of Lemma \ref{prop}.
\end{proof}

 \section{Proof of the main result} 
Presently we establish Theorem \ref{mainth}, which is about the well-posedness of system (\ref{mainprob}), this is done using mainly the result of Lemma \ref{prop}.
\begin{proof}
The proof is carried in two steps. In the first one we construct, starting from a Dirichlet problem, a sequence of functions, $(u^1_\epsilon)_\epsilon$, each of which solves an approximating problem (cf. system \ref{sys11}) which is similar to (\ref{mainprob}) and satisfy the estimate (\ref{propestim}) of Lemma \ref{prop}. In the second subsection we pass to the limit and combine with the a priori estimate (\ref{estim22b}), satisfied by $u^2$, and finally conclude the results of  theorem \ref{mainth}.

\paragraph{\bf Step 1: construction of the approximating solution.} Let $\lbrace s_i\rbrace_{i=1,2}$ be as defined by (\ref{dn}). Denote by $D(s_i, \epsilon_0)\subset\mathbb{R}^2$ the disk centered at $s_i$ with radius $\epsilon_0$, where $\epsilon_0 := \frac{1}{2}\vert s_1 - s_2\vert$. Let $\epsilon$, $0<\epsilon<\epsilon_0$, be a parameter designed to tend towards zero. We introduce a sequence of truncating functions $\chi_{\epsilon, i} \in C^1(\mathbb{R}^2)$, $i =1,2$ defined by:
\begin{equation*}
\chi_{\epsilon, i}(x):=
\begin{cases}
&  0,  \quad  \vert x-s_i \vert\leq \epsilon^2,\\
 &  \exp  \left(\frac{\epsilon -\vert x - s_i \vert}{\epsilon^2 - \vert x - s_i\vert}\right) , \quad  \epsilon^2 <\vert x-s_i\vert <\epsilon,\\
 & 1,  \quad  \epsilon < \vert x-s_i\vert.
 \end{cases}
\end{equation*} 

Define: 
\begin{equation}\label{chi}
\chi_\epsilon(t, x):=
\begin{cases}
 &  \chi_{\epsilon, 1}( x),  \quad  \forall x,  ~ \vert x-s_1\vert\leq \epsilon_0,~ \forall t\in (0,T),\\
 &  \chi_{\epsilon, 2}( x),  \quad  \forall x,  ~  \vert x-s_2\vert\leq \epsilon_0,~ \forall t\in (0,T),\\
 &  1,   \quad  \forall x\in   \mathbb{R}^2 -\mathcal{V}_1\cup \mathcal{V}_2, ~ \forall t\in (0,T),
\end{cases}
\end{equation} 
where we denote $\mathcal{V}_i := D(s_i, \epsilon_0)$. It is easy to see, from the definition (\ref{chi}), that the functions $(\chi_\epsilon)_\epsilon$ are constant with respect to the time variable and satisfy:
\begin{alignat}{2}
1)& ~\vert\chi_\epsilon(t, x)\vert\leq 1,\quad\forall t\in (0,T),\forall x\in\mathbb{R}^2, \forall\epsilon>0.\label{chi1}\\
2)& ~ \epsilon^{1 + \gamma}\Vert\partial_x\chi_\epsilon\Vert_{\infty, (0,T)\times\mathcal{V}_i}\rightarrow 0\quad as ~\epsilon\rightarrow 0, ~ \text{with}
\quad\gamma >0\label{chi2}.\\ 
3)&~ \vert\chi_\epsilon(t, x)\vert\rightarrow 1,\quad as ~\epsilon\rightarrow 0, \forall (t,x)\in(0,T)\times\Gamma_d   \label{chi4}.
\end{alignat}
 Let us cover the set $\overline{\Omega}$ by open sets $U_1$ and $U_2^\epsilon$ such that:
\begin{enumerate}[label=(\roman*)]
\item $\overline{\Omega}\subset  U_1\cup  U^\epsilon_2$, for all $\epsilon>0$\label{i},
\item $U_1 \cap \partial\Omega =  \overline{\Gamma_d}$\label{ii},
\item  $\vert\alpha_i^\epsilon - s_i \vert < \frac{1}{2}\epsilon^2$, for every $i=\overline{1,2}$\label{iii} ,
\end{enumerate}
where we denote $\lbrace \alpha_1^\epsilon, \alpha_2^\epsilon \rbrace := \partial\Omega\cap\partial U_2^\epsilon$ such that $d(\alpha_1^\epsilon, s_1) < d(\alpha_2^\epsilon, s_1)$, for all $\epsilon>0$. The condition (ii) implies that $\partial U_1 \cap\partial\Omega = \lbrace s_1, s_2\rbrace$. Consider, for every $\epsilon>0$, the partition of the unity $(\tilde{\psi}_1^\epsilon, \tilde{\psi}_2^\epsilon)$ subordinate to the cover $(U_1, U_2^\epsilon)$ i.e. 
\begin{enumerate}[label=(\alph*)]
\item $\tilde{\psi}_i^\epsilon \in C_0^\infty( U_i)$, $0<\tilde{\psi}_i^\epsilon <1$,~ for $1\leq i \leq 2$\label{a},
\item $\text{supp}~\tilde{\psi}_1^\epsilon\subset  U_1$, $\text{supp}~\tilde{\psi}_2^\epsilon\subset  U_2^\epsilon$, $\forall\epsilon>0$\label{b},
\item  $\tilde{\psi}_1^\epsilon(x) + \tilde{\psi}_2^\epsilon(x) = 1$ for every $x\in \Omega$\label{c}. 
\end{enumerate}
   Define the functions: 
\begin{equation}
\psi_1^\epsilon(t, x) := \tilde{\psi}_1^\epsilon(x),\quad \psi_2^\epsilon(t,x):= \tilde{\psi}_2^\epsilon(x), \quad\forall t\in (0,T), \forall x\in\Omega .
\end{equation}   
     Thus, the functions $\psi^\epsilon_1, \psi^\epsilon_2 $ are constant with respect to the time variable. One immediate consequence is:  
\begin{equation}\label{convpar}
\Vert\psi_1^\epsilon - 1\Vert_{\infty, (0,T)\times\Gamma_d}\rightarrow 0\quad as~\epsilon\rightarrow 0. 
\end{equation}

Let $G\in H^1\left((0,T)\times\Gamma_d \right)$ and $\tilde{G}\in H^1\left((0,T)\times\partial\Omega\right)$ be such that conditions (\ref{cg1})-(\ref{cg2}) are satisfied. Consider the Dirichlet boundary value problem:
\begin{equation}\label{dirich0}
\left\{
    \begin{array}{ll}
\begin{aligned}
&\partial^2_t u - \Delta u = 0 \quad \text{in}~ (0,T)\times\Omega, \\
&            u = \tilde{G}  \quad \text{on}~ (0,T)\times\partial\Omega,\\      
&           u(0,x) = 0~ , \quad  \partial_tu(0,x) = 0 \quad \text{in}~ \Omega, 
\end{aligned}
   \end{array}
\right.
\end{equation}
the compatibility conditions are well met by $\tilde{G}$ on $\partial\Omega$. Denote $u\in C^0\left( [0,T]; H^1(\Omega)\right)\cap C^1\left( [0,T]; L^2(\Omega)\right)$ to be the unique solution of (\ref{dirich0}). 
It is easy to see that the function $\chi_\epsilon\psi_1^\epsilon u$ satisfy the following system:
\begin{equation}\label{sysys}
\left\{
    \begin{array}{ll}
\begin{aligned}
&\partial^2_t (\chi_\epsilon\psi_1^\epsilon u ) - \Delta(\chi_\epsilon\psi_1^\epsilon u) = f_\epsilon\quad \text{in}~ (0,T)\times\Omega_s, \\
&            \chi_\epsilon\psi_1^\epsilon u = G_\epsilon  \quad \text{on}~  (0,T)\times\Gamma_d,\\      
&          \frac{\partial( \chi_\epsilon\psi_1^\epsilon u)}{\partial \overrightarrow{n}} = 0 \quad  \text{on}~ (0,T)\times\Gamma_n, \\
&     (\chi_\epsilon\psi_1^\epsilon u)(0,x) = 0, \quad \partial_t (\chi_\epsilon\psi_1^\epsilon u)(0,x) = 0~ \quad\text{in}~ \Omega,
\end{aligned}
   \end{array}
\right.
\end{equation}
where $G_\epsilon := \psi_1^\epsilon\chi_\epsilon\tilde{G} $ and $f_\epsilon\in C^1([0,T]; L^2(\Omega))$. One should remark that the homogeneous Neumann condition in (\ref{sysys}) arises from the conditions \ref{ii} and \ref{b}. Consider the function $w_\epsilon$ satisfying the following system:
\begin{equation}\label{sys2}
\left\{
    \begin{array}{ll}
\begin{aligned}
&\partial^2_t w_\epsilon - \Delta w_\epsilon = f - f_\epsilon  \quad\text{in}~ (0,T)\times\Omega_s, \\
&           w_\epsilon  = 0   \quad\text{on}~ (0,T)\times\Gamma_d,\\      
&           \frac{\partial w_\epsilon }{\partial \overrightarrow{n}}  = 0  \quad\text{on}~ (0,T)\times\Gamma_n, \\
&    w_\epsilon(0,x) = \Psi^0(x), \quad \partial_t w_\epsilon(0,x) = \Psi^1(x) \quad\text{in}~ \Omega,
\end{aligned}
   \end{array}
\right.
\end{equation}
where $(\Psi^0, \Psi^1)$ satisfy the weak boundary conditions (\ref{wb}). In virtue of {\cite[Theorem 1, p. 170]{Ibuki}} and its Corollary, see also \cite{Hayashida},  there is a unique  $w_\epsilon \in C^1([0,T]; K(\Omega)) \cap C^2([0,T]; L^2(\Omega))$ satisfying, for all $\epsilon>0$, system (\ref{sys2}). Then it is obvious to remark that the function $u_\epsilon:= \chi_\epsilon\psi_1^\epsilon u + w_\epsilon$ solves the following system:
\begin{equation*}
\left\{
    \begin{array}{ll}
\begin{aligned}
&\partial^2_t u_\epsilon - \Delta u_\epsilon = f  \quad\text{in}~ (0,T)\times\Omega_s, \\
&           u_\epsilon  = G_\epsilon   \quad\text{on}~ (0,T)\times\Gamma_d,\\      
&           \frac{\partial u_\epsilon }{\partial \overrightarrow{n}} = 0  \quad\text{on}~ (0,T)\times\Gamma_n, \\
&           u_\epsilon (0,x) = \Psi^0(x), \quad  \partial_t u_\epsilon (0,x) = \Psi^1(x) \quad\text{in}~ \Omega.
\end{aligned}
   \end{array}
\right.
\end{equation*}
On the other hand, consider the function $u^2$ that solves the system:
\begin{equation}\label{problem2}
\left\{
    \begin{array}{ll}
\begin{aligned}
&\partial^2_t u^2 - \Delta u^2 =  f  \quad\text{in}~ (0,T)\times\Omega,\\
&            u^2 = 0  \quad\text{on}~ (0,T)\times\Gamma_d,\\ 
&         \frac{\partial u^2}{\partial \overrightarrow{n}} = 0  \quad\text{on}~ (0,T)\times\Gamma_n, \\
&          u^2(0, x) = \Psi^0(x), \quad  \partial_t  u^2(0,x) = \Psi^1(x)~ \ \quad\text{in}~\Omega,
\end{aligned}
   \end{array}
\right.
\end{equation}
following again {\cite[Theorem 1, p. 170]{Ibuki}}, see also \cite{Hayashida}, we deduce that (\ref{problem2}) admits a unique solution: $$u^2 \in C^1([0,T]; K(\Omega))\cap C^2([0,T]; L^2(\Omega)), $$ moreover we have, cf. {\cite[Estimate 1.9, p. 170]{Ibuki}}, the following estimate:
\begin{equation}\label{estim22b}
\begin{aligned}
&\Vert u^2\Vert_{C^0([0,T]; H^1(\Omega)) \cap C^1([0,T]; L^2(\Omega))}\\& \leq C (\Vert f\Vert_{L^1(0,T; L^2(\Omega))} +  \Vert\Psi^0\Vert_{H^1(\Omega)} + \Vert\Psi^1\Vert_{L^2(\Omega)}).
\end{aligned}
\end{equation}

\paragraph{\bf Step 2: passing to the limit using Lemma \ref{prop}.} It is immediate to see that $u^1_\epsilon:= u_\epsilon - u^2 \in C^0([0,T]; H^1(\Omega))\cap C^1([0,T]; L^2(\Omega))$ exists and satisfy, for every $\epsilon>0$, the system:
\begin{equation}\label{sys11}
\left\{
    \begin{array}{ll}
\begin{aligned}
&\partial^2_t u_\epsilon^1 - \Delta u_\epsilon^1 =  0\quad\text{in}~ (0,T)\times\Omega;\\
&            u_\epsilon^1 = G_\epsilon \quad\text{on}~ (0,T)\times\Gamma_d;\\ 
&          \frac{\partial u_\epsilon^1}{\partial n}   =  0 \quad\text{on}~(0,T)\times\Gamma_n;     \\
&           u_\epsilon^1(0, x) = 0,\quad   \partial_t   u_\epsilon^1(0,x) = 0~ \text{in}~ \Omega,
\end{aligned}
   \end{array}
\right.
\end{equation}
and thus $G_\epsilon \in\mathcal{G}$ for every $\epsilon>0$, where $\mathcal{G}$ is defined by (\ref{subG}). Moreover, the associated solution $u_\epsilon^1$ of system (\ref{sys11}) satisfy, for every $\epsilon>0$, the assumption of Lemma \ref{prop}. Consequently, if we note $g^1_\epsilon:= \frac{\partial u_\epsilon^1}{\partial n} $, we obtain by applying (\ref{propestim}):
\begin{equation}\label{estim00}
\Vert g^1_\epsilon \Vert_{H^\beta((0,T)\times\Gamma_d)} \leq C_1 \Vert G\chi_\epsilon \psi_1^\epsilon\Vert_{H^1((0,T)\times\Gamma_d)}.
\end{equation}
for all $\epsilon >0$. The $H^1-$norm of $G_\epsilon - G = G\chi_\epsilon \psi_1^\epsilon - G$ satisfy:
\begin{equation}\label{H1norm}
\begin{aligned}
&\Vert G_\epsilon - G\Vert_{H^1((0,T)\times\Gamma_d)}\\  &\leq \Vert G_\epsilon - G\Vert^2_{L^2((0,T)\times\Gamma_d)} + \Vert\partial_t G_\epsilon - \partial_t G\Vert^2_{L^2((0,T)\times\Gamma_d)}\\& +  \Vert\nabla_x G_\epsilon-  \nabla_x  G\Vert^2_{L^2((0,T)\times\Gamma_d)}.
\end{aligned}
\end{equation}
We show that $\Vert G_\epsilon - G\Vert_{H^1((0,T)\times\Gamma_d)} \rightarrow 0$. Since both of $\psi_1^\epsilon$ and $\chi_\epsilon$ are constant with respect to the time variable, and considering (\ref{chi1}), (\ref{chi4}) and (\ref{convpar}), then we can easily see that the limit of the first two terms in the RHS of (\ref{H1norm}) is zero, one can use for instance the dominated convergence theorem. Let us estimate the limit of the term $||\nabla_x G_\epsilon -\nabla_x G||_{L^2((0,T)\times\Gamma_d)}$.
This term involves the $L^2$-norm of the gradient on the surface $ (0,T)\times\partial\Omega $ and therefore we need to apply Estimate (\ref{ruleF}) of Remark \ref{remark} that yield us with expression of $\nabla_x(\psi_1^\epsilon\chi_\epsilon G)$, thus we have: 
\begin{equation}\label{H1norm1}
\begin{aligned}
&\Vert\nabla_x G_\epsilon - \nabla_x G\Vert_{L^2((0,T)\times\Gamma_d)}\\  &\leq  \Vert\chi_\epsilon\psi_1^\epsilon\nabla_x G - \nabla_x G\Vert_{L^2((0,T)\times\Gamma_d)}\\& +  \Vert G\chi_\epsilon\nabla_x\psi_1^\epsilon\Vert_{L^2((0,T)\times\Gamma_d)} +  \Vert G\psi_1^\epsilon\nabla_x\chi_\epsilon\Vert_{L^2((0,T)\times\Gamma_d)}.
\end{aligned}
\end{equation}
Let us note the following three facts:
\begin{enumerate}[label=(\Alph*)]
\item Using Estimates: (\ref{chi1}), (\ref{chi4}) and (\ref{convpar}), we can show easily, using dominated convergence, that:  $$|| \chi_\epsilon\psi_1^\epsilon\nabla_x G - \nabla_x G||_{L^2((0,T)\times\Gamma_d)} \rightarrow 0, \quad\text{as}~\epsilon\rightarrow 0.\label{A} $$

\item Given condition \ref{i}-\ref{iii} and \ref{b}, the support of $\nabla_x\psi^\epsilon_1$ is contained in the region $\cup_i D(s_i, \epsilon^2)$, for every $\epsilon>0$. Moreover, according to its definition (\ref{chi}), the function $\chi_\epsilon$ is zero on this region, thus:
\begin{equation}\label{B}
\begin{aligned}
&\Vert G\chi_\epsilon\nabla_x\psi_1^\epsilon\Vert_{L^2((0,T)\times\Gamma_d)}\\& = \sum_{i=1}^2\Vert G\chi_\epsilon\nabla_x\psi_1^\epsilon\Vert_{L^2\left((0,T)\times \lbrace x\in \Gamma_d,~ \vert x - s_i\vert\leq \epsilon^2 \rbrace\right)}  = 0. 
\end{aligned}
\end{equation}

\item  Using condition (\ref{small}), the properties \ref{a} and (\ref{chi2}), we infer:   
\begin{equation*}
\begin{aligned}
 \Vert G\psi_1^\epsilon\nabla_x\chi_\epsilon\Vert_{L^2((0,T)\times\Gamma_d)} & = \sum_{i=1}^2\Vert G\psi_1^\epsilon\nabla_x\chi_\epsilon\Vert_{L^2\left((0,T)\times \lbrace x\in \Gamma_d,~\epsilon^2 < \vert x -s_i\vert\leq \epsilon \rbrace\right)} \\&  \leq  2\Vert\epsilon^\gamma\Vert_{ L^2\left((0,T)\times\lbrace x\in \Gamma_d,~\epsilon^2 < \vert x -s_i\vert\leq \epsilon \rbrace \right)}    \rightarrow   0, 
\end{aligned}
\end{equation*} 
  we have used the fact that the support of $\partial_x\chi_\epsilon$ on $\Gamma_d$ is contained in $\cup_i\lbrace x\in \Gamma_d,~\epsilon^2 < \vert x -s_i\vert\leq \epsilon \rbrace)$\label{C} .
\end{enumerate}
Consequently, using the facts \ref{A}-\ref{C} we infer that: 
\begin{equation}\label{convpar0}
\Vert\nabla_x (G\psi_1^\epsilon\chi_\epsilon ) - \nabla_x G\Vert_{L^2((0,T)\times\Gamma_d)}\rightarrow 0\quad as~\epsilon\rightarrow 0. 
\end{equation} 
Finally, combining (\ref{H1norm}) and (\ref{convpar0}) one can see that $\Vert G_\epsilon - G\Vert_{H^1((0,T)\times\Gamma_d)}\rightarrow 0$. Consequently, since $\mathcal{A}$ is complete, we deduce from estimate (\ref{estim00}) that there exists $g^1 \in H^\beta((0,T)\times\Gamma_d)$ such that $\Vert g^1_\epsilon - g^1 \Vert_{H^\beta((0,T)\times\Gamma_d)}\rightarrow 0$ and thus, by passing to the limit in estimate (\ref{estim00}), we obtain: 
 \begin{equation}\label{estim000}
\Vert g^1 \Vert_{H^\beta((0,T)\times\Gamma_d)} \leq C \Vert G\Vert_{H^1((0,T)\times\Gamma_d)}.
\end{equation}
Since $G_\epsilon$ and $g^1_\epsilon$ converges strongly in $H^1((0,T)\times\Gamma_d)$ resp. $H^\beta((0,T)\times\Gamma_d)$, then $G_\epsilon \rightarrow G$ a.e. on $(0,T)\times\Gamma_d$  and $g^1_\epsilon\rightarrow 0$ a.e. on $(0,T)\times\Gamma_n$. Using the energy estimate (\ref{energyT}), when associated to the Neumann data $g^1_\epsilon - g^1$ and letting $\epsilon\rightarrow 0$, one can easily see that the limit functions $G$ resp. $g^1$ are the Dirichlet resp. the Neumann data on $ \Gamma_d$ of problem (\ref{mix}) whose solution is $u^1 = \lim u^1_\epsilon$. 
  Combining Estimates (\ref{energy}), (\ref{estim000}) and the continuous embedding of $H^\beta$ in $L^2$ we conclude:
\begin{equation}\label{estim22a}
\begin{aligned}
\Vert u^1\Vert_{C^0(0,T; H^\alpha(\Omega))\cap C^1(0,T; H^{\alpha - 1}(\Omega))} \leq C \Vert G\Vert_{H^1((0,T)\times\Gamma_d)},
\end{aligned}
\end{equation}
 where $u^1$ solves (\ref{mix}) with the function  $G$ as a Dirichlet data on $(0,T)\times\Gamma_d$.
\paragraph{Conclusion.} It is immediate to see that the function $u := u^1 + u^2$ solves uniquely system (\ref{mainprob}), where $u^2$ is the unique solution of problem (\ref{problem2}). Combining estimates (\ref{estim22b}) and (\ref{estim22a}) we infer that:
\begin{equation*}
\begin{aligned}
&\Vert u\Vert_{C^0(0,T; H^\alpha(\Omega))\cap C^1(0,T; H^{\alpha - 1}(\Omega))}\\& = \Vert u^1 + u^2\Vert_{C^0(0,T; H^\alpha(\Omega))\cap C^1(0,T; H^{\alpha - 1}(\Omega))}\\  \leq & C\left( \Vert G\Vert_{H^1\left((0,T)\times\Gamma_d\right)} +  \Vert f\Vert_{L^1\left(0,T; L^2(\Omega)\right)} + \Vert\Psi^0\Vert_{H^1(\Omega)} + \Vert\Psi^1\Vert_{L^2(\Omega)}\right),
\end{aligned}
\end{equation*}
which yields estimate (\ref{propestim22}) that guarantee the stability of the solution. Finally, it is easy to see how the uniqueness of the solution  $u$  follows from this estimate, and thus conclude the proof of Theorem \ref{mainth}.
\end{proof}

\begin{remark}\label{remark}
Since $\tilde{G}\in H^1((0,T)\times\partial\Omega)$ we know, cf. {\cite[Theorem 1, p.40] {Burenkov}},  there exists  $\tilde{G}^n\in C^\infty(\partial\Omega)$ such that $\Vert\tilde{G}^n - \tilde{G}\Vert_{H^1((0,T)\times\partial\Omega)}\rightarrow 0$.
On the other hand, given that $\partial\Omega\in C^\infty $, we know, cf. {\cite[Theorem 2, p.12]{ff}}, that there exists an extension $\tilde{G}_e^n$ of $\tilde{G}^n$ having a sufficient regularity, say $\tilde{G}_e^n\in C^1(\mathbb{R}^2)$ such that: 
\begin{equation}\label{ff}
 \tilde{G}_e^n = \tilde{G}^n \quad\text{on}~ (0,T)\times\partial\Omega
\end{equation}
the rule of derivation of the product on $\Omega$ yields:
\begin{equation}\label{ruleO}
\nabla_x(\tilde{G}_e^n\chi_\epsilon\psi_1^\epsilon) = \chi_\epsilon\psi_1^\epsilon\nabla_x \tilde{G}_e^n + \tilde{G}_e^n\chi_\epsilon\nabla_x\psi_1^\epsilon + \tilde{G}_e^n\psi_1^\epsilon\nabla_x\chi_\epsilon , \quad\text{a.e. in}~(0,T)\times\Omega.
\end{equation}

Given the smoothness of all the involved functions, identity (\ref{ruleO}) holds true when passing to the boundary points $\Gamma_d$. Thus, taking into account (\ref{ff}), we have:
\begin{equation}\label{ruleB}
\nabla_x(G^n\chi_\epsilon\psi_1^\epsilon) = \chi_\epsilon\psi_1^\epsilon\nabla_x G^n + G\chi_\epsilon\nabla_x\psi_1^\epsilon + G^n\psi_1^\epsilon\nabla_x\chi_\epsilon , \quad\text{a.e. in}~(0,T)\times\Gamma_d.
\end{equation}

Using the density assumption, we pass to the limit $n\rightarrow\infty$ by using the dominated convergence theorem to deduce that:
\begin{equation}\label{gloi}
\begin{aligned}
&\Vert\nabla_x(G^n\chi_\epsilon\psi_1^\epsilon) - \left(\chi_\epsilon\psi_1^\epsilon\nabla_x G +  G\chi_\epsilon\nabla_x\psi_1^\epsilon  G\psi_1^\epsilon\nabla_x\chi_\epsilon \right)\Vert_{H^1((0,T)\times\partial\Omega}\\& \leq  \Vert\chi_\epsilon\psi_1^\epsilon\nabla_x G^n - \chi_\epsilon\psi_1^\epsilon\nabla_x G\Vert +  \Vert G^n\chi_\epsilon\nabla_x\psi_1^\epsilon- G\chi_\epsilon\nabla_x\psi_1^\epsilon \Vert\\& + \Vert G^n\psi_1^\epsilon\nabla_x\chi_\epsilon  - G\psi_1^\epsilon\nabla_x\chi_\epsilon  \Vert\\&\rightarrow 0.
\end{aligned}
\end{equation}

On one hand using the density assumption, and on the other hand converting the normed convergences in (\ref{gloi}) into point wise convergences, we get using (\ref{ruleB}): 
\begin{equation}\label{ruleF}
\begin{aligned}
 \nabla_x G_\epsilon(t,x) &=  \nabla_x(G\chi_\epsilon\psi_1^\epsilon)(t,x)\\ &=   \lim_{n\rightarrow\infty} \nabla_x(G^n\chi_\epsilon\psi_1^\epsilon)(t,x)\\& = \chi_\epsilon\psi_1^\epsilon\nabla_x G + G\chi_\epsilon\nabla_x\psi_1^\epsilon + G\psi_1^\epsilon\nabla_x\chi_\epsilon , \quad\text{a.e. in}~(0,T)\times\Gamma_d.
 \end{aligned}
\end{equation}
\end{remark}

\begin{remark}\label{31}
One should notice that the derivative $\partial_x \psi^\epsilon_1$ can become infinite when $\epsilon\rightarrow 0$, the point behind the domination condition (\ref{small}) and condition (iii) is to annihilate the effect of this possible blow-up.


\noindent - Finally, we emphasize that the results stated in \cite{Lions} and which are used in this paper, were established within the framework of the Laplace operator. That said, they remain valid in case of a strongly and uniformly elliptic operator $ \mathcal{L} $ in (\ref{hyperbol}).
\end{remark}


\end{document}